\documentclass{amsart}
\usepackage{palatino,hhline, bbm}
\usepackage{arydshln} 
\usepackage{amsthm, amsmath}
\usepackage{amssymb} 
\usepackage{amscd}
\usepackage{mathrsfs}
\renewcommand{\mathcal}{\mathscr}
\usepackage[hmargin=3.1cm, vmargin=3cm]{geometry}
\usepackage[breaklinks,bookmarksopen,bookmarksnumbered]{hyperref}
\usepackage[all]{xy}

\theoremstyle{plain}
\newtheorem{thm}{Theorem}
\newtheorem{lem}{Lemma}
\newtheorem{prop}{Proposition}
\newtheorem*{cor}{Corollary}

\theoremstyle{remark}

\newcommand\pr{\noindent\textit{Proof} : }

\newcommand\rond{\kern 1pt{\scriptstyle\circ}\kern 1pt}

\def\lr#1{\langle {#1} \rangle}

\newcommand\codim{\operatorname{codim}}

\newcommand\Pic{\operatorname{Pic}}

\newcommand\pp{_{}^{\scriptscriptstyle \bullet}}
\newcommand{\mo}{\smallsetminus}
\newcommand\dep{\operatorname{depth}}

\newcommand\Z{\mathbb{Z}}

\newcommand\C{\mathbb{C}}
\renewcommand\P{\mathbb{P}}

\renewcommand\O{\mathcal{O}}

\newcommand\iso{\vbox{\hbox to .8cm{\hfill{$\scriptstyle\sim$}\hfill}
\nointerlineskip\hbox to .8cm{{\hfill$\longrightarrow $\hfill}} }}
\newcommand\bir{\vbox{\hbox to .8cm{\hfill{$\scriptstyle\sim$}\hfill}
\nointerlineskip\hbox to .8cm{{\hfill$\dasharrow $\hfill}} }}

\begin{document}
\title{An introduction to Ulrich bundles}
\author[Arnaud Beauville]{Arnaud Beauville}
\address{Universit\'e C\^ote d'Azur\\
CNRS -- Laboratoire J.-A. Dieudonn\'e\\
Parc Valrose\\
F-06108 Nice cedex 2, France}
\email{arnaud.beauville@unice.fr}
 
\begin{abstract}
After recalling the definition and basic properties of Ulrich bundles, we focus on the existence problem: does any smooth projective variety carry a Ulrich bundle?  We show that the Serre construction provides a positive answer on certain surfaces and threefolds.
\end{abstract}
\maketitle 
\section*{Introduction}
Ulrich bundles have had three lives. They first appeared in the 80's
 in commutative algebra, starting with \cite{U}, under various names (linear maximal, or maximally generated, Cohen-Macaulay modules). They entered the realm of algebraic geometry twenty years later, with the beautiful paper \cite{ES}. But it is only in the recent years that they have received wide  attention among algebraic geometers.

  In these notes we will first explain the definition and basic properties of Ulrich bundles, and motivate their study by the case of hypersurfaces. We will then focus on the fundamental and  intriguing question of the subject, which goes back to \cite{U}: does any smooth projective variety carry a Ulrich bundle?  We will show that the classical Serre  construction of rank 2 bundles give interesting examples for some surfaces and threefolds.

\medskip	
We will work over the complex numbers, though most  results are valid over any algebraically closed field. We will denote the ambient projective space by $\P$ when we do not need to specify its dimension. 
\bigskip	
\section{Ulrich bundle on hypersurfaces}

We start with a  classical problem: given a smooth hypersurface $X\subset \P$, can we write the equation of $X$ as a determinant  of linear forms? The (positive) answer for a smooth cubic surface  goes back at least to  1863 \cite{S}. But  determinantal 
hypersurfaces have singularities in codimension 3, so a smooth hypersurface of dimension $\geq 3$ cannot be written as a determinant.

Let us settle therefore for a weaker property: we ask whether $X$ can be defined \emph{set-theoretically} by  a linear determinant. This can be expressed nicely in terms of vector bundles on $X$:

\begin{prop}
Let $X$ be a smooth hypersurface of degree $d$ in $\P$, given by an equation $F=0$, and let $r$ be an integer $\geq 1$. The following conditions are equivalent:

$1)$ $F^r=\det(L_{ij})$, where $(L_{ij})$ is a $rd\times rd$ matrix of linear forms on $\P$;

$2)$ There exists a rank $r$ vector bundle $E$ on $X$ and an exact sequence $0\rightarrow \O_{\P}(-1)^{rd}\xrightarrow {\ L\ }\O_{\P}^{rd}\rightarrow E\rightarrow 0$.
\end{prop}
\pr If $2)$ holds, the determinant of $L$ vanishes exactly along $X$, hence is proportional to some power of $F$, thus to $F^r$ for degree reasons.

Assume that $1)$ holds. We claim that the cokernel $E$ of the homomorphism $L:\O_{\P}(-1)^{rd}\rightarrow \O_{\P}^{rd}$ is a vector bundle on $X$. Let $x\in X$; the depth of  the stalk $E_x$ as a $\O_{\P,x}$-module or as a $\O_{X,x}$-module is the same, and is equal by the Auslander-Buchsbaum formula to\[\dep(E_x)=\dim(\O_{\P,x})-\mathrm{proj.dim}_{\O_{\P,x}}(E_x)=\dim(\O_{\P,x})-1=\dim(\O_{X,x})\ ,\]
hence  $E_x$ is a free $\O_{X,x}$-module, so  $E$ is locally free on $X$.

The first Chern class of $E$ as a $\O_{\P}$-module is  $ \ -rdc_1(\O_{\P}(-1))= r[X]$, hence $E$ has rank $r$.  \qed

\medskip	
The vector bundles $E$ which appear in $2)$ will be our main object of study. They can be characterized in a variety of ways:
\begin{prop}
Let $X\subset\P  $ be a smooth hypersurface of degree $d$, $E$ a rank $r$ vector bundle on $X$. The following conditions are equivalent:

$1)$ There exists an exact sequence $0\rightarrow \O_{\P}^{rd}(-1)\xrightarrow {\ L\ }\O_{\P}^{rd}\rightarrow E\rightarrow 0$.

$2)$ The cohomology $H\pp(X,E(-p))$ vanishes for $1\leq p\leq \dim(X)$.

$3)$ If $\pi :X\rightarrow H$ is the projection from a point of $\P\mo X$ to a hyperplane, the vector bundle $\pi _*E$ is trivial (hence isomorphic to $\O_H^{rd}$). 
\end{prop}

It turns out that this Proposition is a particular case of a more general result, which makes sense on \emph{any} smooth projective variety:

\begin{thm}[\cite{ES}]
Let $X\subset \P$ be a smooth variety, and let $E$ be a vector bundle on $X$. The following conditions are equivalent:

$1)$ There exists a linear resolution 
\[0\rightarrow L_c\rightarrow L_{c-1}\rightarrow \ldots \rightarrow L_0\rightarrow E\rightarrow 0\]with $c=\codim (X,\P)$ and  $L_i=\O_{\P}(-i)^{b_i}$.

$2)$ The cohomology $H\pp(X,E(-p))$ vanishes for $1\leq p\leq \dim(X)$.

$3)$ If $\pi :X\rightarrow \P^{\dim(X)}$ is a \emph{finite} linear projection, 
 the vector bundle $\pi _*E$ is trivial.

The vector bundle $E$ is a \emph{Ulrich bundle} if it satisfies these equivalent conditions. 
\end{thm}

\pr We first observe that the Theorem holds for $X=\P^n$; that is, a vector bundle $E$ on $\P^n$ with $H\pp(E(-p))=0$ for $1\leq p\leq n$  is trivial. Indeed this implies $H^p(\P^n,E(-p))=0$ for $p>0$, which means that $E$ is $0$-\emph{regular} in the sense of Mumford \cite[Lecture 14]{M}; this implies  in particular that $E$ is globally generated, and satisfies $H^p(\P^n,E)=0$ for $p>0$. The Hilbert polynomial $\chi (E(t))$ vanishes for $t=-1,\ldots ,-n$, and its leading coefficient is $\frac{r}{n!} $; therefore $\chi (E(t))=\frac{r}{n!} (t+1)\ldots (t+n)$, hence $h^0(E)=r$. Thus we have a surjective homomorphism $\O_{\P^n}^r\rightarrow E$, which is necessarily an isomorphism.

\medskip	
Let us treat the general case. We put $n:=\dim (X)$. If 3) holds, we have  $H\pp(X,E(-i))=\allowbreak H\pp(\P^{n}, \pi _*E(-i))=0$ for $1\leq i\leq n$, hence 2). Conversely, if 2) holds, the vector bundle $F:=\pi _*E$ on $\P^n$ satisfies  $H\pp(\P^n,F(-i))=0$ for $1\leq i\leq n$, hence is trivial by the case already treated. 

\medskip	
The implication $1) \ \Rightarrow\ 2)$ follows from the vanishing of $H\pp(\P,\O_{\P}(-i))$ for $1\leq i\leq \dim(\P)$. Assume that $2)$ holds. We will define by induction a sequence of $0$-regular sheaves $K_i$ on $\P$, for $0\leq i\leq c$, such that:

$a)$ $K_{0}=E$;

$b)$ $K_{i+1}(-1)$ is the kernel of the evaluation map 
$H^0(\P,K_i)\otimes \O_{\P}\rightarrow K_i$;

$c)\ H\pp(K_i(-j))=0$ for $1\leq j\leq n+i$.

Suppose the $K_i$ are defined for $0\leq i\leq p$; we define $K_{p+1}$ by $b)$.  From the exact sequence $0\rightarrow K_{p+1}(-1)\rightarrow H^0(\P,K_p)\otimes \O_{\P}\rightarrow K_p\rightarrow 0$ and $c)$ we get $H\pp(K_{p+1}(-j))=0$ for $1\leq j\leq n+p+1$, 
and also $H^q(K_{p+1}(-q))=H^{q-1}(K_{p}(-(q-1)))=0$. Thus $K_{p+1}$ is $0$-regular and satifies $b)$ and $c)$.

Put $L_i:=H^0(\P,K_i)\otimes \O_{\P}(-i) $; the exact sequences $0\rightarrow K_{i+1}(-i-1)\rightarrow L_i\rightarrow K_i(-i)\rightarrow 0$  give a long exact sequence
\[0\rightarrow K_{c}(-c)\rightarrow L_{c-1}\rightarrow L_{c-1} \ldots L_0\rightarrow E\rightarrow 0\ .\]
Now $c)$ means that $K_c$ is a Ulrich bundle on $\P$, hence trivial by the first part of the proof. \qed

\bigskip	
\noindent\emph{Remarks}$.-$ 1) Condition $2)$ depends only on the line bundle $\O_X(1)$ (and not on the space of sections giving the embedding in $\P$). 
We will sometimes say that $E$ is a Ulrich bundle for $(X,\O_X(1))$.

 2) A direct sum
of vector bundles is a Ulrich bundle iff each summand is a Ulrich bundle. Thus we need only to classify  indecomposable Ulrich bundles.

 \bigskip	
Using the theorem, we can now answer our question on representation of hypersurfaces by linear determinants. Already in the case of quadrics, the result and the proof are far from trivial:

\begin{prop}\label{q}
Let $Q\subset \P^{n+1}$ be a smooth quadric. If $n$ is odd (resp. even) there is exactly one (resp. two) indecomposable Ulrich bundles on $Q$, of rank $2^{[\frac{n-1}{2} ]}$. 

As a consequence, we can write
$(X_1^2+\ldots +X_n^2)^r$ as a determinant of linear forms if and only if $r$ is a multiple of $2^{[\frac{n-3}{2}]}$. 
\end{prop}

\smallskip	
\pr The projection $\pi :Q\rightarrow \P^n$ from a point  $p\in\P^{n+1}\mo X $ realizes $Q$ as a double covering of $\P^n$. We choose the coordinates so that $Q$ is given by an equation $X^2_{n+1}-q(X_0,\ldots ,X_n)=0$, and $\pi (X_0,\ldots ,X_{n+1})=(X_0,\ldots ,X_n)$. Thus $\pi $ is branched along the quadric $q=0$ in $\P^n$. This implies  $\pi _*\O_Q=\O_{\P^n}\oplus \O_{\P^n}(-1)$, with the multiplication $\O_{\P^n}(-1)\otimes \O_{\P^n}(-1)\rightarrow \O_{\P^n}$ given by $q$.

The data of a vector bundle $E$ on $Q$   is equivalent to the data of the vector bundle $\pi _*E$ on $\P^n$, together with an algebra homomorphism 
$\varphi :\pi _*\O_Q\rightarrow \mathcal{E}nd(\pi _*E)$. If $E$ is a rank $r$ Ulrich bundle, we have $\pi _*E\cong \O_{\P^n}^{2r}$, and $\varphi $ is given by a matrix $A\in \mathbf{M}_{2r}(\C)\otimes_{\C} H^0(\P^n,\O_{\P^n}(1))$, that is, a matrix $A= (A_{ij}(X))$ whose entries are linear forms in $X=(X_0,\ldots ,X_n)$, satisfying $A^2=q(X)\cdot I_{2r}$.   

Equivalently, we can view $A$ as a $\C$-linear map $\tilde{A}: \C^{n+1}\rightarrow \mathbf{M}_{2r}(\C)$, satisfying $\tilde{A}(v)^2=q(v)\cdot I_{2r} $ for all $v$ in $\C^{n+1}$. This is  the same as a $\C$-algebra  homomorphism $C_{n+1}\rightarrow \mathbf{M}_{2r}(\C)$, where $C_{n+1}$ is the Clifford algebra $C(\C^{n+1},q)$. 
If $n$ is odd,   the algebra $C_{n+1}$ is simple, hence has a unique simple module, of dimension $2^{\frac{n+1}{2} }$, which corresponds to an indecomposable Ulrich bundle of rank $r=2^{\frac{n-1}{2} }$. 
If $n$ is even, $C_{n+1}$ is a product of two simple algebras of dimension $2^{n}$, hence has two simple modules of dimension $2^{\frac{n}{2} }$, corresponding to two non-isomorphic indecomposable Ulrich bundle of rank $2^{\frac{n-2}{2} }$. \qed

\medskip	
\noindent\emph{Remark}$.-$ The bundles $\Sigma $ (resp.\ $\Sigma _+$ and $\Sigma _-$) appear in \cite{K}, with a stronger result: the vector bundles $\lr{\Sigma ,\O_Q,\ldots ,\O_Q(n-1)}$ for $n$ odd, $\lr{\Sigma_{-} ,\Sigma _{+},\O_Q,\ldots ,\O_Q(n-1)}$ for $n$ 
even form a semi-orthogonal decomposition of the (bounded) derived category of coherent sheaves on $Q$.

\bigskip	
Using the notions of generalized Clifford algebra associated to a higher degree form and of matrix factorization, Herzog, Ulrich and Backelin have obtained
a vast generalization of Proposition \ref{q}:
\begin{thm}[{\cite{HUB}}]
Any smooth complete intersection $X\subset\P$ carries a Ulrich bundle.
\end{thm}

\bigskip	
\section{Some general properties}
We first list a few direct consequences of Theorem 1.

\noindent (2.1)  Let $E$ be a rank $r$ Ulrich vector bundle on $X\subset \P$. Then $H^i(X,E(j))=0$ for all $j$ and $0<i<\dim(X)$ and  $h^0(E)=rd$ (this follows  from the condition 3) of the Theorem).

\smallskip	
\noindent (2.2) A Ulrich bundle is semi-stable; if if it is not stable, it is an extension of Ulrich bundles of smaller rank 
(this follows   again from condition 3).

\smallskip	
\noindent (2.3)  The Ulrich bundles on a curve $C\subset\P$ are the bundles $E(1)$, where $E$ is a vector bundle with vanishing cohomology.

\smallskip	
\noindent (2.4) If $E$ is a Ulrich bundle on $X\subset \P$ and $Y$ is a hyperplane section of $X$, $E_{|Y}$  is a Ulrich bundle on $Y$.
This follows from the exact sequence $0\rightarrow E(-1)\rightarrow E\rightarrow E_{|Y}\rightarrow 0$.

\smallskip	
\noindent (2.5) Let $E$ and $F$ be  Ulrich bundles for $(X,\O_X(1))$ and $(Y,\O_Y(1))$; put $n:=\dim(X)$. Then $E\boxtimes F(n)$ is a Ulrich bundle for $(X\times Y, \O_X(1)\boxtimes  \O_Y(1))$.
 
\hskip0.5cm  Indeed we have $H\pp(X\times Y,E(-p)\boxtimes F(n-p))=H\pp(X,E(-p))\otimes H\pp(Y,F(n-p))$. The first factor is zero for $1\leq p\leq n$ and the second one for $n+1\leq p\leq n+\dim(Y)$.

\smallskip	
\noindent (2.6) Let $\pi :X\rightarrow Y$ be a finite surjective morphism,  
$L$ a very ample line bundle on $Y$, $E$ a vector bundle on $X$. Then $E$ is a Ulrich bundle for $(X,\pi ^*L)$ if and only if $\pi _*E$ is a Ulrich bundle for $(Y,L)$.

\hskip0.5cm This  follows  from the isomorphism $H\pp(Y, \pi _*E\otimes L^{-k})\cong H\pp(X, E\otimes \pi ^*L^{-k})$.

\begin{prop}\label{ver}
$(\P^n,\O_{\P^n}(d))$ admits a Ulrich bundle of rank $n!$.
\end{prop}
\pr Consider the quotient map $\pi :(\P^1)^n\rightarrow \mathrm{Sym}^n\P^1=\P^n$; the pull-back $\pi ^*\O_{\P^n}(1)$ is \break $\O_{\P^1}(1)\boxtimes\ldots \boxtimes \O_{\P^1}(1)$. By (2.5) the line bundle $L=  \O_{\P^1}(d-1)\boxtimes\ldots \boxtimes \O_{\P^1}(nd-1)$ is a Ulrich line bundle for $((\P^1)^n, \pi ^*\O_{\P^n}(d))$,  hence $\pi _*L$ is a Ulrich bundle for  $(\P^n,\O_{\P^n}(d))$ by  (2.6). \qed

\medskip	
In \cite{ES} the authors use a much more sophisticated method to construct a \emph{homogeneous} Ulrich bundle for $(\P^n,\O_{\P^n}(d))$, of rank $d^{\binom{n}{2}}$. 
\begin{cor}[\cite{ES}]
Let $X\subset \P$ be a smooth variety of dimension $n$, carrying  a Ulrich bundle of rank $r$. Then $(X,\O_X(d))$ carries a Ulrich bundle of rank $rn!$.  
\end{cor}
\pr Let $E$ be  a Ulrich bundle of rank $r$ for $X$, $\pi :X\rightarrow \P^n$  a finite projection, and  $F$  a Ulrich bundle of rank $n!$ for $(\P^n,\O_{\P^n}(d))$. Then $E\otimes \pi ^*F$ is a Ulrich bundle of rank $rn!$ for $(X,\O_X(d))$.\qed

\bigskip	
\section{Ulrich line bundles}
Ulrich line bundles are rather exceptional; for instance, they cannot exist on a variety $X\subset \P$ with $\Pic(X)=\Z\,\O_X(1)$, except when $X$ is a linear subspace. Apart from the case of curves and  projective spaces, there are a few situations where they do exist:

\begin{prop}
$1)$ Let $S\subset \P$ be a del Pezzo surface (embedded by its anticanonical system). Let $L$ be a line bundle on $S$ with $L^2=-2$, $L\cdot K_S=0$. Then $L(1)$ is a Ulrich line bundle on $S$. We can take for instance $L=\O_S(\ell-\ell')$, where $\ell$ and $\ell'$ are two disjoint lines.

\smallskip	
$2)$ Let $X\subset \P$ be a scroll over a curve $B$: that is, there exists a fibration $p:X\rightarrow B$ whose fibers are linear subspaces in $\P$.   Let $M$ be a  line bundle on $B$ with $H\pp(B,M)=0$. Then $p^*M\otimes \O_X(1)$ is a Ulrich line bundle on $X$. 
\end{prop}

\pr 1) We have $\chi (L)=0$ by Riemann-Roch, $H^0(S,L)=0$ since $L\cdot K_S=0$, and $H^2(S,L)=0$ since $(K _S\otimes L^{-1}\cdot K_S^{-1})<0$, hence $H\pp(S,L)=0$. Then
$H\pp(S,L^{-1})=0$ for the same reason, hence 
 $H\pp(S,L(-1))=0$ by Serre duality.

\smallskip	
2) We have $Rp_*(p^*M\otimes \O_X(-k))=M\otimes Rp_*\O_X(-k)$. The object $Rp_*\O_X(-k)$ of $D(B)$ vanishes for $k=1,\ldots ,\dim(X)-1$; for $k=0$ we get $H\pp(X, p^*M)\cong H\pp(B,M)=0$.\qed
 
 \medskip	
 Another interesting case is that of Enriques surfaces $S\subset\P$ \cite{BN}. Suppose $S$ has no $(-2)$-curve. If we can express $\O_S(1)$ as $\O_S(A-B)$, where $A,B$ are divisors with $A^2=B^2=-2$, then $\O_S(A)\otimes \O_S(1)$ is a Ulrich line bundle, because $\O_S(A)$ and $\O_S(A)\otimes \O_S(-1)\cong \O_S(B)$ have vanishing cohomology. Borisov and Nuer conjecture that this is always possible, and prove it for the Fano model $S\subset\P^5$.

\bigskip
\section{Rank 2 Ulrich bundles on surfaces}
Let $S\subset\P$ be a smooth projective surface. 
To find Ulrich bundles on $S$, a natural attack is to look for bundles $E$ such that $E(-1)$ and $E(-2)$ are Serre dual, so that the vanishing of $H\pp(E(-1))$ implies that of $H\pp(E(-2))$. If $E$ has rank 2, this is achieved by imposing  $\det E=K_S(3)$, which means that $E$ is \emph{special} in the sense of \cite{ES}: we say that a rank 2 Ulrich bundle $E$ on a  smooth projective variety $X$ of dimension $n$ is \emph{special} if $\det E=K_X(n+1)$. The existence of such a bundle implies that the Chow form of $X$ admits a ``B\'ezout expression" as a Pfaffian (\emph{loc. cit.}). 
\begin{prop}\label{kod0}
Minimal surfaces of Kodaira dimension $0$ admit a special rank $2$ Ulrich bundle, except perhaps some special K3 surfaces.
\end{prop}
\pr Recall that there are four classes of such surfaces, namely K3, Enriques, abelian and bielliptic surfaces.
The case of K3 surfaces is treated in \cite{AFO}, using the Lazarsfeld-Mukai construction. For each polarization type, the special surfaces for which the method does not apply form a strict subvariety  of the moduli space. 

\medskip	
We will give a uniform proof for the three remaining cases\footnote{The case of abelian surfaces was treated in \cite{Bab}, and that of Enriques surfaces in \cite{Bpq}
and \cite{C}.}, using  the Serre construction  in the following form (see for instance \cite[Theorem 5.1.1]{HL}). Let $S\subset \P$ be a smooth surface; using Remark 1 (\S 1) we can assume that $S$ is embedded by the complete linear system $|\O_S(1)|$.
Recall that a finite set of points $Z\subset S$ has the \emph{Cayley-Bacharach property} if for every $p\in Z$, any hyperplane containing $Z\mo\{p\} $ contains $Z$. If this is the case, there exists a rank 2 vector bundle $E$ and an extension
\[0\rightarrow K_S\rightarrow E\rightarrow \mathcal{I}_Z(1)\rightarrow 0\ . \leqno{(*)}\]
Our aim is to show that a certain twist of $E(1)$ is a Ulrich bundle. We will use the following easy lemma\,:

\begin{lem}\label{lem}
Let $S\subset \P$ be a smooth surface, and let $E$ be a rank $2$ vector bundle on $S$ with $\det E=K_S(1)$, $h^0(E)=\chi (E)=0$. Then $E(1)$ is a special Ulrich bundle.
\end{lem}
\pr We have $K_S\otimes E^*\cong E(-1)$, hence $h^2(E)=h^0(E(-1))=0$. Since $\chi (E)=0$ this implies  $H\pp(E)=0$, then  $H\pp(E(-1))=H\pp(K_S\otimes E^*)=0$.\qed

\medskip	
As a consequence, we have the following result, due to Casnati \cite{C}, and improving \cite{Bpq}:
\begin{prop}\label{pq0}
Let  $S\subset \P^n$ be a surface with $q=p_g=0$ and $H^1(S,\O_S(1))=0$. Then $S$ admits a  rank $2$ special Ulrich bundle. 
\end{prop}
\pr We choose for $Z$ a set of $n+2$ general points of $S$. They form a projective frame of $\P^n$, so $Z$ satisfies the Cayley-Bacharach property.  We have $\chi (\mathcal{I}_Z(1))=\chi (\O_S(1))-n-2=-1$, hence $\chi (E)=0$. We have $H^0(\mathcal{I}_Z(1))=0$ since $Z$ is general, hence $E(1)$ is a Ulrich bundle by  Lemma~\ref{lem}.\qed

\medskip	
Let us come back to Proposition \ref{kod0}. Proposition \ref{pq0} implies the case of Enriques surfaces. Let $S\subset \P^n$ be an abelian or bielliptic surface. We choose a smooth hyperplane section $C=S\cap H$ of $S$, and a subset of $n+1$ general points of $C$. Then any subset of $n$ points of $Z$ span $H$, hence $Z$ satisfies the Cayley-Bacharach property. Note that $g(C)=n+2$ by Riemann-Roch.

Let $\eta $ be an element of order 2 of $\Pic^{\mathrm{o}}(S)$, with $\eta \neq \O_S,K_S$; we will apply Lemma \ref{lem} to $E\otimes \eta $. Since $\chi (\O_S)=0$ and $\chi (\mathcal{I}_Z(1))=\chi (\O_S(1))-n-1=0$, we have $\chi (E\otimes \eta )=\chi (E)=0$. Let us twist the exact sequence $(*)$ by $\eta $; we have $H\pp(K_S\otimes \eta )=H\pp(\eta )=0$, hence it suffices to prove $H^0(\mathcal{I}_Z\otimes \eta (1))=0$. Since the restriction map $H^0(S,\eta (1))\rightarrow H^0(C,\eta (1)_{|C})$ is bijective, it suffices to show $H^0(C,\eta (1)_{|C}(-Z))=0$. We have $\eta (1)_{|C}=K_C\otimes (\eta \otimes K_S^{-1})_{|C}$; 
since $\eta \otimes K_S^{-1}$ is nontrivial, so is its restriction to $C$, hence   $h^0(C, \eta (1)_{|C})=g(C)-1=n+1$ and  $h^0(C,\eta (1)_{|C}(-Z))=0$ because $Z$ is general. Thus by Lemma~\ref{lem} $E\otimes \eta (1)$ is a special Ulrich bundle.

\medskip	
\noindent\emph{Remark}$.-$  It is quite possible that every rational surface carries a rank 2 Ulrich bundle.
On the other hand,  as observed in \cite[\S 6]{ES}, a general surface $S\subset \P^3$ of degree $\geq 16$ does \emph{not} carry a rank 2 Ulrich bundle.
Indeed the existence of such a bundle on a surface $S$ with $\Pic(S)=\Z$ implies that the equation of $S$ can be written as a Pfaffian of linear forms \cite[Corollary 2.4]{Bdet}, and this is not possible for a general $S\subset \P^3$ of degree $\geq 16$ \cite[Proposition 7.6]{Bdet}.

\bigskip	
\section{Fano threefolds of index $2$}

Recall that a (smooth) Fano threefold $X$ has even index  if there exists a line bundle $L$ on $X$ with $L^2\cong K_X^{-1}$. If $L$ is very ample, it  embeds $X$ as a subvariety of degree $d$ of $\P^{d+1}$, with $d:=(L^3)$. There exist 7 families of such threefolds, with $3\leq d\leq 8$ \cite{F}; for $d=8$ we have $X=\P^3$ and $L=\O_{\P^3}(2)$, otherwise $X$ has index 2, that is, $L$ is not divisible  in $\Pic(X)$.

\begin{prop}
Every Fano threefold  $X\subset \P^{d+1}$ of even index   admits a special rank $2$ Ulrich bundle.
\end{prop}

\pr Recall that a complete linear system of degree $n\geq 3$ on an elliptic curve $\Gamma $  defines an embedding $i:\Gamma \hookrightarrow \P^{n-1}$; the curve $i(\Gamma )$ is projectively normal, and is called a normal elliptic curve. We will prove in the Lemma below that the threefold $X$ contains
 a normal elliptic curve $\Gamma \subset \P^{d+1}$ (of degree $d+2$). Since $\omega _{\Gamma }\cong \O_{\Gamma }$, we have $\det N_{\Gamma /X}\cong \O_{\Gamma }(2)$. Since $H^2(X,\O_X(2))=0$, there exists
by the Serre construction (see e.g.\ \cite{A})  a rank 2 vector bundle $E$ on $X$ and an exact sequence
\[0\rightarrow \O_X\rightarrow E\rightarrow \mathcal{I}_{\Gamma }(2)\rightarrow 0\ .\leqno{(*)}\]
Let us show that $E$ is a Ulrich bundle. Since $\det E=\O_X(2)$,  we have $E(-2)\cong K_X\otimes E(-2)^*$ and $E(-3)\cong K_X\otimes E(-1)^*$. Therefore using Serre duality it suffices to prove $H\pp(E(-1))=0$ and $H^i(E(-2))=0$ for $i=0,1$. But this follows from the equalities
\[H\pp(\O_X(-1))=H\pp(\mathcal{I}_{\Gamma }(1))=0\quad\mbox{and}\quad H^i(\O_X(-2))=H^i(\mathcal{I}_{\Gamma })=0\ \mbox{ for } i=0,1.\eqno{\qed}\]

\begin{lem}
$X$ contains a normal elliptic curve $\Gamma \subset \P^{d+1}$.
\end{lem}
\pr Assume first $d\leq 7$. Let $H$ be a general hyperplane in $\P^{d+1}$; the surface $S:=X\cap H$ is a Del Pezzo surface, obtained by blowing up $9-d$ points of $\P^2$ in general position. Consider the linear system of quartic curves in $\P^2$ passing doubly through two of these points, and simply through the others. It is easy to see that it is base point free, and therefore contains a smooth curve, which gives a normal elliptic curve  $\Gamma _0\subset H$ of degree $d+1$.

Now let $p$ be a point of $\Gamma _0$, and let $L$ be a line contained in $X$ but not in $H$, and passing through $p$.  We want to show that the curve $C=\Gamma _0\cup L$ can be deformed to a smooth curve $\Gamma $ in $X$, which must be a normal elliptic curve of degree $d+2$. According to \cite[Theorem 4.1]{HH}, it suffices to prove  $H^1(L, N_{L/X})=0$ and
$H^1(\Gamma_0,  N_{\Gamma _0/X}(-p))=0$. The normal bundle $N_{L/X}$ is equal to $\O_L^2$ or $\O_L(1)\oplus \O_L(-1)$ \cite[Lemma 3.3.4]{F}, hence its $H^1$ vanishes.
The normal bundle $N_{\Gamma _0/X}$ is an extension of $(N_{S/X})_{|\Gamma _0}=\O_{\Gamma _0}(1)$ by $N_{\Gamma _0/S}$; since both line bundles have degree $d+1\geq 4$, we have  $H^1(\Gamma_0,  N_{\Gamma _0/X}(-p)) =0$, hence our assertion. 
 
 Now we consider the case $d=8$, so  $X=\nu (\P^3)\subset \P^9$, where $\nu $ is the Veronese embedding given by the linear system $|\O_{\P^3}(2)|$. 
 Let $\Gamma $ be a  smooth elliptic curve of degree 5 in $\P^3$.
 We observe that  a smooth quadric or a quadratic cone cannot contain an elliptic curve of degree 5, so $\Gamma $ is not contained in a quadric; therefore the restriction map
 $\nu^*: H^0(\P^3,\O_{\P^3}(2))\rightarrow H^0(\Gamma ,\O_{\Gamma }(2))$ is bijective, and $\nu (\Gamma )$ is a normal elliptic curve in $\P^9$.\qed

\medskip	
\noindent\emph{Remarks}$.-$ 1) The existence of a Ulrich bundle (possibly of higher rank) could also be deduced case by case from the explicit description of these threefolds given in \cite{F}.
 
 \smallskip	
2) In the particular case of cubic threefolds, the Ulrich bundles given by the Proposition have already been constructed and studied in \cite{MT}. 
 They have a beautiful geometry, see \cite{MT} and also \cite{D} and \cite{Bcub}. They are also considered for $d=4$ or $5$ in \cite{A-C}.

\smallskip
3) Conversely, any  special Ulrich bundle $E$ of rank 2 on $X$ is of the above form. Indeed, let $s$ be a general section of $E$; the zero locus $\Gamma $ of $s$ is a smooth curve, and we have an exact sequence
\[0\rightarrow \O_X\xrightarrow{\ s\ } E\rightarrow \mathcal{I}_{\Gamma }(2)\rightarrow 0\ .\]
Since $E$ is Ulrich, we have $h^1(\mathcal{I}_{\Gamma })=0$ and $h^2(\mathcal{I}_{\Gamma })=1$, hence $h^0(\O_{\Gamma })=h^1(\O_{\Gamma })=1$, so that $\Gamma $ is an elliptic curve. Moreover we have $h^0(\mathcal{I}_{\Gamma }(1))=h^1(\mathcal{I}_{\Gamma }(1))=0$, hence $\Gamma $ is a normal elliptic curve in $\P^{d+1}$.

\smallskip	
4)  If $3\leq d\leq 5$, we have $\Pic(X)=\Z$, so any rank 2 Ulrich bundle is special and stable, because there are no Ulrich line bundles on $X$ (see (2.2)). It is easy to see that this is still the case for $d=7$ or $8$. However for 
$d=6$ there exist non-special rank 2 Ulrich bundles, and non-stable special rank 2 Ulrich bundles.

\smallskip	
5) The same method produces a rank 2 Ulrich bundle on the Fano fourfolds of index 3 $X\subset \P^{d+2}$ containing a del Pezzo surface of degree $d+2$. However these Fano fourfolds are not general. For instance,  the cubic fourfolds containing a degree 5 del Pezzo surface are exactly the pfaffian cubics, that is, those which can be defined by the pfaffian of a $6\times 6$ skew-symmetric matrix of linear forms; they form a hypersurface in the moduli space of cubic fourfolds \cite[Proposition 9.2]{Bdet}. 

\bigskip	
We now consider the moduli space   of rank 2 special Ulrich bundles  on $X$;  it is  an open subset of the moduli space of   semi-stable rank 2 bundles.

\begin{prop}
The moduli space $\mathcal{M}$ of rank 2 special Ulrich bundles  on $X$ is smooth of dimension $5$.
\end{prop}
 \pr Let $\mathcal{H}$ be the Hilbert scheme of normal elliptic curves in $\P^{d+1}$ contained in $X$. Since $H^1(X,\O_X(2))=0$, a curve $\Gamma $ in $\mathcal{H}$ determines the vector bundle $E$ in $\mathcal{M}$ up to isomorphism \cite{A}. This defines a map $p: \mathcal{H}\rightarrow \mathcal{M}$ which is surjective (Remark 3 above); the fiber $p^{-1}(E)$ is canonically identified with an open subset of $\P(H^0(X,E))$, so $p$ is smooth of relative dimension $2d-1$. 

 Let $\Gamma \in\mathcal{H}$, and let $N$ be the normal bundle of $\Gamma $ in $X$. From the exact sequence $(*)$ we get an isomorphism $N\cong E_{|\Gamma }$. Therefore $\deg N=2d+4$. We will  prove that $H^1(\Gamma ,N)=0$; this implies that the Hilbert scheme $\mathcal{H}$ is smooth of dimension $2d+4$, hence that $\mathcal{M}$  is smooth of dimension $5$.
 
 The Ulrich bundle $E$ admits a presentation
 \[\O_{\P}(-1)^q\xrightarrow{\ M\ } \O_{\P}^p\rightarrow E\rightarrow 0\]
 where $M$ is a matrix of linear forms. Restricting to $\Gamma $, dualizing and taking cohomology, we get an exact sequence
 \[0\rightarrow H^0(\Gamma , N^*)\rightarrow H^0(\Gamma ,\O_{\Gamma })^p\xrightarrow{\  ^t\!M\ } H^0(\Gamma ,\O_{\Gamma }(1))^q\ .\]
 Let $v\in H^0(\Gamma ,N^*)$; we have $^t\!Mv=0$ in $H^0(\Gamma ,\O_{\Gamma }(1))^q\cong H^0(\P^{d+1} ,\O_{\P }(1))^q$. Since $M$ is generically surjective, $^t\!M$ is generically injective, so this implies $v=0$. Therefore $H^0(\Gamma ,N^*)=0$, and $H^1(\Gamma ,N)=0$ by Serre duality.\qed
 
 \medskip	
\noindent \emph{Examples}$.-$ 1) Let $X$ be a cubic threefold, and let $J^2X$ be the translate of the intermediate Jacobian of $X$ which parametrizes 1-cycles of degree 2. 
The map $\mathcal{M}\rightarrow J^2X$ which associates to a rank 2 Ulrich bundle $E$ the class $c_2(E(-1))$ is an open embedding; its image is the complement  of the divisor in $J^2X$ parametrizing sum of two lines \cite[Corollary 6.4]{Bcub}.

2) Consider $X=\P^3$, embedded in $\P^9$ by the Veronese embedding. Any rank 2 Ulrich bundle on $X$  is the  kernel of a homomorphism $T_{\P^3}(1)\rightarrow \O_{\P^3}(3)$ given by a contact form $\alpha  \in H^0(\P^3, \Omega ^1(2))$ \cite[Proposition 5.11]{ES}. Thus the moduli space $\mathcal{M}$ is the open subset of contact forms in $\P(H^0(\P^3, \Omega ^1(2)))$, or equivalently the space of non-degenerate skew-symmetric forms on $\C^4$ up to a scalar. 

\medskip	
\noindent\emph{Remark}$.-$ Assume $3\leq d\leq 5$. Let $S$ be a very general surface in the linear system $|\O_X(2)|$, so that $S$ is a K3 surface with $\Pic(S)=\Z$. Let $\mathcal{M}_S$ be the moduli space of stable vector bundles on $S$ with $c_1=0$, $c_2=4$; it is an open subset of the holomorphic symplectic  10-dimensional manifold constructed by O'Grady \cite{O}. The proof of  \cite[Proposition 8.4]{Bcub} applies identically to show   that \emph{the map $E\mapsto E(-1)_{|S}$ induces an isomorphism of $\mathcal{M}$ onto a Lagrangian subvariety of} $\mathcal{M}_S$. However for $d=4$ or $5$ these subvarieties do not seem to be part of a Lagrangian fibration.

\end{document}